\documentclass{amsart}

\usepackage{epsfig}

\newtheorem{theorem}{Theorem}[section]
\newtheorem{lemma}[theorem]{Lemma}
\newtheorem{corollary}[theorem]{Corollary}
\newtheorem{proposition}[theorem]{Proposition}

\newtheorem{definition}[theorem]{Definition}

\newtheorem{remark}[theorem]{Remark}

\newcommand{\re}{\mathrm{reg}}
\newcommand{\al}{\mathrm{al}}
\newcommand{\lp}{\ell^\prime}
\newcommand{\wh}{\widehat}
\newcommand{\wt}{\widetilde}

\newcommand{\rh}{\mathrm{rh}}
\newcommand{\inv}{\mathrm{inv}}
\newcommand{\diag}{\mathrm{diag}}
\newcommand{\supp}{\mathrm{supp}}
\newcommand{\nbc}{\mathrm{NBC}}
\newcommand{\CM}{\nbc(w)}
\newcommand{\ot}{\leftarrow}
\newcommand{\bg}{\mathrm{bg}}
\newcommand{\Sn}{\mathfrak{S}_n}
\newcommand{\Sm}{\mathfrak{S}_m}

\newcommand{\R}{\mathbb{R}}
\newcommand{\A}{\mathcal{A}}
\newcommand{\SP}{\mathrm{span}}
\newcommand{\out}{\mathcal{E}}

\begin{document}
\author{Axel Hultman}
\address{Department of Mathematics, KTH-Royal Institute of Technology, 
  SE-100 44, Stockholm, Sweden.}
\title{Inversion arrangements and Bruhat intervals}

\begin{abstract}
Let $W$ be a finite reflection group. For a given $w\in W$, the following
assertion may or may not be satisfied:
\vspace{3pt}
\begin{itemize}
\item[($*$)]The principal Bruhat order ideal of $w$ contains as many
  elements as there are regions in the inversion hyperplane arrangement of $w$.
\end{itemize}
\vspace{3pt}
We present a type independent combinatorial criterion which
characterises the elements $w\in W$ that 
satisfy ($*$). A couple of immediate consequences are derived:
\begin{itemize}
\item[(1)] The criterion only involves the order ideal of $w$ as an
  abstract poset. In this sense, ($*$) is a poset-theoretic property. 
\item[(2)] For $W$ of type $A$, another characterisation of ($*$), in
  terms of pattern avoidance, was previously given in collaboration
  with Linusson, Shareshian and Sj\"ostrand. We obtain
  a short and simple proof of that result.
\item[(3)] If $W$ is a Weyl group and the Schubert variety
  indexed by $w\in W$ is rationally smooth, then $w$ satisfies ($*$).
\end{itemize}
\end{abstract}

\maketitle

\section{Introduction}
Let $n$ be a positive integer. Given indices $1\le i < j \le n$,
define a hyperplane
\[
H_{i,j} = \{(x_1, \ldots, x_n)\in \R^n\mid x_i=x_j\}.
\]
The arrangement of all such hyperplanes
\[
\A_n = \{H_{i,j}\mid 1\le i<j\le n\}
\]
is known as the {\em braid arrangement}. The orthogonal reflections in
the hyperplanes $H_{i,j}$ generate a finite reflection group
isomorphic to the symmetric group $\Sn$; a natural isomorphism is
given by associating a reflection through $H_{i,j}$ with the transposition
$(i, j) \in \Sn$.

Given a permutation $w\in \Sn$, we define its {\em inversion
  arrangement} as the following subarrangement of $\A_n$:
\[
\A_w = \{H_{i,j}\mid  1\le i<j\le n, \, w(i)>w(j)\}.
\]
In particular, $\A_{w_0} = \A_n$, where $w_0\in \Sn$ is the reverse permutation
$i\mapsto n+1-i$.

The inversion arrangement cuts the ambient space into a set $\re(w)$ of
{\em regions}, a region being a connected component of the complement $\R^n \setminus \cup\A_w$. 

Let $[\cdot,\cdot]$ denote closed intervals in the Bruhat order on
$\Sn$ (the definition of which is recalled in Section
\ref{se:prel}). Postnikov \cite{postnikov} discovered a numerical relationship between $\re(w)$ and the Bruhat order ideal $[e,w]$, where $e\in \Sn$ is the
identity permutation. When $w$ is a Grassmannian permutation, he
proved that the sets are equinumerous; both are in 1-1 correspondence with
certain cells in a CW decomposition of the totally nonnegative
Grassmannian. For arbitrary $w$, he conjectured the following results that were
subsequently proven in \cite{HLSS}:
\begin{itemize}
\item[(A)] For all $w\in \Sn$, $\#\re(w) \le \#[e,w]$.
\item[(B)] Equality holds in (A) if and only if $w$ avoids the patterns $4231$,
  $35142$, $42513$ and $351624$.
\end{itemize}
The reader who is not familiar with the terminology employed in (B)
finds an explanation in Section \ref{se:newproof}. 

We have just defined $\A_w$ using $\Sn$-specific language. It
is, however, completely natural to replace $\Sn$ by an arbitrary finite
reflection group $W$ and consider $\A_w$, $\re(w)$ and $[e,w]$ for any $w\in
W$; see Section \ref{se:prel} for details of the definitions. In fact,
it was not (A) but the following result which was established in \cite{HLSS}:

\begin{itemize}
\item[(A$^\prime$)] Given a finite reflection group $W$ and any $w\in
  W$, $\#\re(w) \le \#[e,w]$.  
\end{itemize}

This generalises (A),\footnote{An explanation of
  the implication (A$^\prime$) $\Rightarrow$ (A) can be found in
  \cite{HLSS}.} but notice that there is no statement 
(B$^\prime$). Indeed, the problem of how to characterise those $w\in W$
for which equality holds in (A$^\prime$) was posed as
\cite[Open problem 10.3]{HLSS}. Such a characterisation is the main result of the present
paper. The precise assertion is stated in Theorem 
\ref{th:main}. It essentially says that equality holds in (A$^\prime$)
if and only if the following property is satisfied for every $u\le w$:
among all paths of shortest length from $u$ to $w$ in the Cayley graph
of $W$ (with edges generated by reflections), there is one which visits
vertices in order of increasing Coxeter length. 

A number of consequences are derived from the main
result:

First, we conclude that the characterising property is
poset-theoretic. That is, whether or not equality holds in (A$^\prime$)
can be determined by merely looking at $[e,w]$ as an abstract poset.

Second, we give a new proof of the difficult direction of (B). In
\cite{HLSS}, (A$^\prime$) was proven by exhibiting an
injective map $\phi$ from (essentially) $\re(w)$ to $[e,w]$. Thus,
proving (B) amounts to characterising surjectivity of $\phi$ in terms
of pattern avoidance when $W=\Sn$. That surjectivity implies the
appropriate pattern avoidance is a reasonably straightforward consequence
of the construction of $\phi$; 
see \cite[Section 4]{HLSS}. Contrastingly, the proof of the converse
statement given in \cite[Section 5]{HLSS} is a direct, fairly
involved, counting argument which does not use $\phi$ at all. In light of
our Theorem \ref{th:main}, surjectivity of $\phi$ can now, however, be
related to pattern avoidance in a rather straightforward way.

Third, when $W$ is a Weyl group, each element $w\in W$ corresponds to a
Schubert variety $X(w)$. We derive from Theorem 
\ref{th:main} that equality holds in (A$^\prime$)
whenever $X(w)$ is rationally smooth. To this end, we establish a
variation of the classical Carrell-Peterson criteria for rational
smoothness which should be of independent interest. It is to be noted
that Oh and Yoo \cite{OY} recently derived a stronger $q$-analogue
of equality in (A$^\prime$) for rationally smooth $X(w)$.

Here is an outline of the structure of the remainder of the paper. In
the next section we agree on basic notation and concepts related to
reflection groups. In
particular, the definition of the map $\phi$ is recalled from
\cite{HLSS}. In Section \ref{se:main}, we 
establish our main result. The new proof of (B) is described in
Section \ref{se:newproof} before we conclude in Section
\ref{se:ratsmooth} with the connection to rationally smooth Schubert
varieties.
 
\section{Reflection groups and inversion arrangements} \label{se:prel}
In this section, we recall some properties of finite reflection
groups. The reader looking for more information should consult
\cite{BB} or \cite{humphreys}. We also review parts of \cite{HLSS}
that are needed for subsequent sections. 

A finite reflection group $W$ is the same as a finite {\em Coxeter
  group}. It is generated by a set $S$ of {\em simple reflections}
  subject to relations of the form $s^2=e$ for all $s\in S$ and
  $(ss^\prime)^{m(s,s^\prime)}=e$ for suitable $m(s,s^\prime)$. Here,
  $e\in W$ is the identity element. 

For $w\in W$, the {\em Coxeter length} $\ell(w)$ is the smallest $k$
such that $w=s_1\cdots s_k$ for some $s_i\in S$. The expression
$s_1\cdots s_k$ is then called {\em reduced}. 

The set $T$ of {\em reflections} consists of all conjugates of simple
reflections, i.e.\ $T=\{wsw^{-1}\mid w\in W\}$. The {\em absolute
  length} $\lp(w)$ is the smallest $k$ such that $t_1\cdots t_k=w$ for
some $t_i\in T$. 

Choose a root system $\Phi\subset \R^n$ for $W$ with set of positive roots
$\Phi^+$. In an incarnation of $W$ as a group generated by orthogonal
reflections in Euclidean space, the
positive roots are in one-to-one correspondence with the reflections
of $W$; the reflecting hyperplane fixed by a reflection is the orthogonal
complement of the corresponding root.

When $W$ is a symmetric group $\Sn$, so that $T$ is the set of
transpositions, it is well known that $\lp(w) = n-c(w)$, where $c(w)$
is the number of cycles in the disjoint cycle decomposition of
$w$. This fact is generalised by the following fundamental
result of Carter which connects the absolute length function with the
underlying geometry.

\begin{theorem}[Carter \cite{carter}] \label{th:carter}
Let $W$ be a finite reflection group. Given $w\in W$, the following
assertions hold.
\begin{itemize}
\item[(a)] The codimension of the fixed point space of $w$ equals
  $\lp(w)$.
\item[(b)] Given reflections $t_1, \ldots, t_m\in T$, we have
  $\lp(t_1\cdots t_m)=m$ if and only if the corresponding roots
  $\alpha_{t_1}, \ldots, \alpha_{t_m}\in \Phi^+$ are linearly independent.
\end{itemize}
\end{theorem}

\begin{remark}\label{re:carter}
{\em A useful consequence is that if there are two
minimal factorisations into reflections $t_1\cdots t_m = r_1\cdots r_m = w$,
$\lp(w)=m$, then we must have $\SP\{\alpha_{t_1}, \ldots, \alpha_{t_m}\} =
\SP\{\alpha_{r_1}, \ldots, \alpha_{r_m}\}$ since both sides of the
equality sign coincide with the orthogonal complement of the fixed
point space of $w$.} 
\end{remark}

The {\em Bruhat graph} $\bg(W)$ is the Cayley graph of $W$ with edges directed towards greater Coxeter length. That is, the vertex set is $W$ and we have directed edges $x\to tx$ for $x\in W$, $t\in T$, whenever $\ell(x)<\ell(tx)$.

\begin{figure}[t]
\begin{center}
\epsfig{height=6cm, file=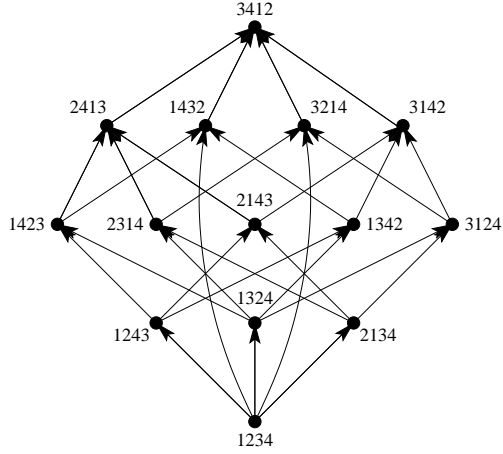}
\caption{The Bruhat graph associated with the permutation $3412\in
  \mathfrak{S}_4$ (one line notation). Disregarding the two curved
  edges yields the Hasse diagram of the Bruhat interval $[e=1234,
  3412]$.} \label{fi:bruhat_graph}
\end{center}
\end{figure}

Taking transitive closure of $\bg(W)$ yields the {\em Bruhat order} on
$W$. In other words, $u\le w$ if and only if $u\to \cdots \to w$. The
subgraph of $\bg(W)$ which is induced by the principal order ideal
$[e,w] = \{u\in W\mid u\le w\}$ is denoted by $\bg(w)$. We refer to
$\bg(w)$, too, as a Bruhat graph. An example can be found in Figure \ref{fi:bruhat_graph}.

Let $\al(u,w)$ denote the distance from $u$ to $w$ in $\bg(w)$
(equivalently, in $\bg(W)$) in the directed, graph-theoretic
sense. Thus, $\al(u,w)$ is finite precisely when $u\le w$. Clearly,
$\al(u,w)\ge \lp(uw^{-1})$ in general, since the right hand side can
be thought of as the distance from $u$ to $w$ in $\bg(W)$ if we
disregard directions of edges. 

A convenient characterisation of the Bruhat order can be given in terms of reduced expressions:

\begin{proposition}
Choose a reduced expression $s_1\cdots s_k$ for $w\in W$. Then, $u\le w$ if and only if $u=s_1\cdots \wh{s_{i_1}}\cdots \wh{s_{i_m}}\cdots s_k$ for some $1\le i_1 < \cdots < i_m\le k$, where a hat denotes omission of an element.
\end{proposition}

The equivalence of these two appearances of the Bruhat order can be derived from the following fundamental fact.

\begin{proposition}[Strong exchange property]
If $u\to w$ and $s_1\cdots s_k$ is any expression for $w\in W$, then
$u = s_1\cdots \wh{s_i} \cdots s_k$ for some $i\in [k]=\{1, \ldots, k\}$.
\end{proposition}

For the remainder of this section, $s_1\cdots s_k$ is a fixed reduced expression for some $w\in W$, where $W$ is a finite Coxeter group. The {\em inversions} of $w$ are the reflections of the form $t_i = s_1s_2\cdots s_{i-1}s_is_{i-1}\cdots s_2s_1$, $i\in [k]$. The
set $\inv(w)$ of inversions of $w$ is independent of the choice of
reduced expression.

Let $\alpha_i \in \Phi^+$ be the root corresponding to $t_i$, and
denote by $H_i = \alpha_i^\perp$ the associated hyperplane. The {\em
  inversion arrangement} of $w$ is
\[
\A_w = \{H_1, \ldots, H_k\}.
\]
The connected components of the complement of $\cup \A_w$ are called {\em
  regions} of $\A_w$. The set of such regions is denoted by $\re(w)$.

At the heart of \cite{HLSS} one finds the construction of an injective map $\re(w) \to [e,w]$. (More accurately, the domain of
the map is not $\re(w)$, but a set which is equinumerous with
$\re(w)$.) We shall study this map further in the present paper, so we
review it here. For convenience, we deviate slightly from the presentation
in \cite{HLSS}, but the formulations are equivalent via
standard facts from matroid theory. 

It is convenient to order positive roots that correspond to inversions
of $w$ with respect to the indices. For example, $\{\alpha_{i_1} <
\cdots < \alpha_{i_m}\}$ indicates the set $\{\alpha_{i_1}, \ldots,
\alpha_{i_m}\}$ under the assumption $1\le i_1<\cdots<i_m\le k$. 

A {\em circuit} is a minimal linearly dependent set $X = \{\alpha_{i_1}< \cdots < \alpha_{i_m}\}\subseteq \Phi^+$ of positive roots corresponding to inversions of $w$ in the manner described above. If $X$ is a circuit, $\{\alpha_{i_1}< \cdots < \alpha_{i_{m-1}}\}$ is a {\em broken circuit}.\footnote{Note that a broken circuit is a circuit missing its {\em largest} element. This convention is backwards compared to common matroid terminology but convenient for our purposes.} If
$Y\subseteq \{\alpha_1, \ldots, \alpha_k\}$ does 
not have a subset which is a broken circuit, say that $Y$ is an {\em
  NBC set}. We denote the family of NBC sets by $\nbc(w)$, although it
of course depends not only on $w$ but also on the choice of reduced
expression $s_1\cdots s_k$. The point
is that $\#\re(w) = \#\nbc(w)$. This well known fact follows for
instance by combining 
two different interpretations of the characteristic polynomial of
$\A_w$ evaluated at $-1$. The $\re(w)$ part of the story is due to Zaslavsky
\cite{zaslavsky} whereas the $\nbc(w)$ connection in this generality
was presented by Rota \cite{rota}. 

\begin{definition}
Construct a map $\phi:\nbc(w) \to [e,w]$ by $\{\alpha_{i_1}< \cdots <
\alpha_{i_m}\} \mapsto t_{i_1}\cdots t_{i_m}w$.
\end{definition} 

Proving statement (A$^\prime$), it was shown in \cite{HLSS} that
$\phi$ always is well defined and injective. 

\section{A surjectivity characterisation}\label{se:main}
Maintain the notation of the previous section. Thus, we keep fixed a
finite Coxeter group $W$, an element $w\in W$ with a reduced
expression $s_1\cdots s_k$ and corresponding inversions $t_i$ with
their associated positive roots $\alpha_i$, $i\in [k]$.

In this section, we determine when the map $\phi$ is
surjective. The image of $\phi$ is dependent on the choice of reduced
expression for $w$, but the cardinality of the image is not, since it coincides
with $\#\re(w)$. Thus, whether or not $\phi$ is surjective depends
solely on the element $w$.

The next lemma is the main source from which this paper flows.

\begin{lemma}\label{le:preimage}
Assume $\al(u,w) = \lp(uw^{-1})$ for all $u\le w$. For fixed $u\le w$,
let $m = \al(u,w)$ and choose a lexicographically
maximal sequence $(i_m, \ldots, i_1)$ such that $u = s_1\cdots \wh{s_{i_1}}\cdots \wh{s_{i_m}}
\cdots s_k$.\footnote{By the strong exchange property, such a sequence
  exists.} Then, $\{\alpha_{i_1} < \cdots < \alpha_{i_m}\} \in \nbc(w)$.
\begin{proof}
Suppose $u$ is such that the indices $1\le i_1 < \cdots < i_m\le k$ yield a
counterexample with $m$ 
minimal. This minimality implies that if $(j_{m-1}, \ldots, j_1)$ is
lexicographically maximal 
such that $s_1\cdots \wh{s_{j_1}}\cdots \wh{s_{j_{m-1}}} \cdots s_k =
s_1\cdots \wh{s_{i_1}}\cdots \wh{s_{i_{m-1}}} \cdots s_k$, then $\{\alpha_{j_1}
< \cdots < \alpha_{j_{m-1}}\}\in \nbc(w)$. 

If $j_{m-1}=i_m$, then $uw^{-1} = t_{j_1}\cdots t_{j_{m-2}}$ and,
consequently, $\lp(uw^{-1})\le m-2$ which is a contradiction. Thus,
$j_{m-1}\neq i_m$. 

Define $V = \SP\{\alpha_{i_1}, \ldots, \alpha_{i_m}\}$. By Carter's
result (Theorem \ref{th:carter}), $\dim V = m$. Let  
\[
n = \max\{i\in[k]\mid\alpha_i\in V\}.
\]
We claim that $n>i_m$. If $j_{m-1}>i_m$, this is immediate
since $\alpha_{j_{m-1}}\in V$ by Remark \ref{re:carter}. If, on the other hand, $j_{m-1}<i_m$, we have $i_x=j_x$ for all $x\in [m-1]$ by maximality of
$(i_m, \ldots, i_1)$. Any broken circuit which is a subset of
$\{\alpha_{i_1},\ldots,\alpha_{i_m}\}$ therefore contains
$\alpha_{i_m}$. By assumption, such a broken circuit exists, and the
claim is established.

Having concluded $n>i_m$, observe
$uw^{-1}t_nw = s_1\cdots \wh{s_{i_1}}\cdots \wh{s_{i_m}} \cdots
\wh{s_n}\cdots s_k \le w$. Again by Carter's result, $\lp(uw^{-1}t_n)
\le m$. Multiplication by a reflection changes the
absolute length by exactly one, so we conclude $\lp(uw^{-1}t_n)
= m-1$. Thus, $uw^{-1}t_n = t_{a_1} \cdots t_{a_{m-1}}$ for some NBC
set $\{\alpha_{a_1}<\cdots<\alpha_{a_{m-1}}\}\subset V$. By Remark
\ref{re:carter}, $V = \SP \{\alpha_{a_1},\ldots,\alpha_{a_{m-1}},
\alpha_n\}$. Thus, $a_{m-1}<n$ and the
fact that $uw^{-1} = t_{a_1} \cdots t_{a_{m-1}}t_n$ therefore
contradicts maximality of the sequence $(i_m, \ldots, i_1)$. 
\end{proof}
\end{lemma}

The desired characterisation is now within reach. For symmetric
groups, it was established in \cite[Theorem 6.3]{HLSS}. The general
case answers \cite[Open problem 10.3]{HLSS}. 

\begin{theorem}\label{th:main}
The map $\phi:\CM\to [e,w]$ is surjective, hence bijective, if and
only if $\al(u,w)=\lp(uw^{-1})$ for all $u\le w$.
\begin{proof}
The only if part is a direct consequence of the following ``going-down
property'' of $\phi$ (\cite[Proposition 6.2]{HLSS}): if $\phi$ is surjective, then the NBC
set $\phi^{-1}(u) = \{\alpha_{i_1}<\cdots <\alpha_{i_m}\}$,
$m=\lp(uw^{-1})$, corresponds to reflections $t_{i_1}, \ldots,
t_{i_m}\in T$ such that $t_{i_{j-1}}\cdots t_{i_m}w \to t_{i_j}\cdots
  t_{i_m}w$ for all $j$. This immediately implies $\al(u,w)=m$.

Under the assumption $\al(u,w)=\lp(u,w)$ for all $u\le w$, Lemma \ref{le:preimage} provides a preimage $\phi^{-1}(v)$ for any $v\le w$, thereby
establishing the if direction. 
\end{proof}
\end{theorem}

When looking for a shortest path, in the undirected sense, from $u$ to
$w$ in the Bruhat graph, we {\em a priori} have to
consider all of $\bg(W)$. Fortunately, the situation is simpler than
that; the next lemma implies, in particular, that an undirected path
from $u$ to $w$ of length $\lp(uw^{-1})$ can be found inside $\bg(w)$
if $u\le w$.

\begin{lemma}\label{le:path}
Given any $u,w\in W$, there exists an element
$v\le u,w$ such that $\al(v,w) + \al(v,u) = \lp(uw^{-1})$.
\begin{proof}
The Bruhat subgraph induced by a coset corresponding to a reflection
subgroup $D = \langle t_1,t_2\rangle \subseteq W$, where $t_1, t_2\in T$, is
isomorphic to the Bruhat graph of the dihedral Coxeter group which is
isomorphic to $D$ \cite{dyer}. The simple structure of such Bruhat
graphs shows that whenever $x\to y \ot z$, there
exists some $y^\prime$ with $x\ot y^\prime \to z$. This implies that,
in the Bruhat graph $\bg(W)$,
among all (not necessarily directed) paths from $u$ to $w$ of fixed
length $l$, those that are 
minimal with respect to the sum of the Coxeter lengths of the
vertices are of the form $u=x_0 \ot x_1 \ot \cdots \ot x_k \to x_{k+1} \to
\cdots \to x_l = w$ for some $0\le k\le l$. If we let $l =
\lp(uw^{-1})$, $v=x_k$ is an element with the prescribed properties.
\end{proof}
\end{lemma}

As an example, one readily verifies that the directed distance from
any vertex to the top element always coincides with the undirected
distance in Figure \ref{fi:bruhat_graph}. Thus, $\phi$ is surjective
when $w=3412\in \mathfrak{S}_4$. This, of course, is also immediate from
the pattern avoidance condition in statement (B).

An interesting consequence is that $\#\re(w)=\#[e,w]$ is a combinatorial 
property of the poset $[e,w]$. In the symmetric group setting, this
was established in \cite[Corollary 6.4]{HLSS}. 
  
\begin{theorem}
If $w, w^\prime \in W$ satisfy $\#\re(w)=\#[e,w]$ and
$\#\re(w^\prime)<\#[e,w^\prime]$, then $[e,w] \not \cong [e,w^\prime]$
as posets. 
\begin{proof}
Dyer \cite{dyer} has shown that the Bruhat graph $\bg(w)$ is
determined by the combinatorial structure of $[e,w]$. By Lemma
\ref{le:path}, it is therefore possible to determine from the poset
structure of 
$[e,w]$ whether or not $\al(u,w) = \lp(uw^{-1})$ for all $u\le
w$. Invoking Theorem \ref{th:main}, that is sufficient for deciding
whether $\phi$ is surjective.
\end{proof}
\end{theorem}

\section{The symmetric group case revisited}\label{se:newproof}

We interpret composition of permutations from left to right. That is, $uw(i) =
w(u(i))$ for $u,w\in \Sn$, $i\in [n]$.\footnote{When $W=\Sn$, this
  makes our concept of inversions (defined in Section
  \ref{se:prel}) coincide with 
  that which is standard for permutations. Composing from right to
  left would require minor adjustments in the proofs, but not
  in the results.} 

For permutations $p\in \Sm$ and $w\in \Sn$, say that $w$
{\em contains the pattern} $p$ if there exist indices
$1\le i_1<\cdots<i_m\le n$ such that for all $1\le j < k \le m$,
$p(j)<p(k)$ if and only if $w(i_j)<w(i_k)$. If $w$ does not contain
the pattern $p$, it {\em avoids} $p$. 

If $w\in \Sn$ avoids the patterns 4231,
35142, 42513 and 351624, then $\#[e,w] = \#\re(w)$. This is
the difficult direction of statement (B); the fairly involved proof
given in \cite{HLSS} is based on 
deriving a common recurrence relation for $\#[e,w]$ and $\#\re(w)$ and
does not use any properties of the map $\phi$. Finding a direct proof of
surjectivity of $\phi$ was formulated as \cite[Open problem
  10.1]{HLSS}. The purpose of this section is to derive such a proof
from the results of the previous section.

\begin{figure}[t]
\begin{center}
\epsfig{height=2.7cm, file=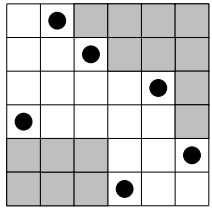}
\end{center}
\caption{A set of dots forming the diagram of the permutation $235164\in
  \mathfrak{S}_6$. The right hull $\rh(235164)$ consists of the
  non-shaded squares.} \label{fi:righthull}
\end{figure}

We shall use a characterisation of the permutations
that avoid the four patterns which is due to Sj\"ostrand
\cite{sjostrand}. To this end, define the {\em diagram} of a
permutation $w\in \Sn$ as the set $\diag(w) = \{(i,w(i))\mid i\in
[n]\} \subset [n]^2$. We think of it as a set of dots on an $n\times
n$ chessboard with matrix conventions for row and column indices, so
that, for instance, $(1,1)$ is the upper left square. 
\begin{definition}
Given $w\in \Sn$, the {\em right hull}  $\rh(w)$ is the subset of
$[n]^2$ which consists of those
$(i,j)$ such that each of the rectangles $\{(x,y)\mid x\le i, \,
y\ge j\}$ and $\{(x,y)\mid x\ge i, \,y\le j\}$ has nonempty
intersection with $\diag(w)$. 
\end{definition}
These concepts are illustrated in Figure \ref{fi:righthull}.

For $w\in \Sn$ and $i,j\in [n]$, let 
\[
w[i,j] = \#\{x\in [n]\mid x \le i, \, w(x)\ge j\}.
\]
The Bruhat order on a symmetric group has the
following convenient characterisation, a proof of which can be found
e.g.\ in \cite{BB}:

\begin{proposition}\label{pr:standard}
For $u,w\in \Sn$, $u\le w$ if and only if $u[i,j]\le w[i,j]$
for all $i,j\in [n]$.
\end{proposition}

Taking into account that $180^\circ$ diagram rotation yields a Bruhat order automorphism, Proposition \ref{pr:standard} makes it clear that $u\le w$ implies
$\diag(u)\subseteq \rh(w)$. We are interested in the permutations $w$
that satisfy the converse. 

\begin{theorem}[Sj\"ostrand \cite{sjostrand}] \label{th:sjostrand}
For $w\in \Sn$, the following are equivalent:
\begin{itemize}
\item $w$ has the {\em right hull property}, meaning $[e,w] =
  \{u\in \Sn\mid \diag(u) \subseteq \rh(w)\}$. 
\item $w$ avoids $4231$, $35142$, $42513$ and $351624$.
\end{itemize}
\end{theorem}

This section is motivated by the desire to find a simple new proof of
(B), so since we are going to use Theorem \ref{th:sjostrand} in that
process, it is relevant to note that Sj\"ostrand's proof (in part
based on ideas of Gasharov and Reiner \cite{GR}) is both elegant and
conceptual.

In light of Theorem \ref{th:sjostrand} and our main result, the if part of
(B) now is equivalent to the following statement:

\begin{lemma}\label{le:RHP}
 If $w\in \Sn$ has the right hull
property, then $\al(u,w)=\lp(uw^{-1})$ for all $u\le w$.
\begin{proof}
Assume $w$ has the right hull property and pick $u < w$. To argue by
induction, it suffices to find a transposition $t\in T$ such that
$u\to tu\le w$ and $\lp(uw^{-1}) = \lp(tuw^{-1}) + 1$. 

Choose a nontrivial cycle $c$ in the disjoint cycle decomposition 
of $uw^{-1}$. Then, $cw < w$ because every dot in the diagram of
$cw$ also appears either in the diagram of $w$ or in that of $u$, both
of which are contained in $\rh(w)$.  

Let $\supp (c) = \{i_1<\cdots<i_m\} \subseteq [n]$ denote the set of non-fixed
elements of $c$. Defining
\[
\mathfrak{S}_c = \{x\in \Sn\mid x(i) = w(i) \text{ for all } i \not
\in \supp (c)\}, 
\]
we thus have $w, cw\in \mathfrak{S}_c$. A natural bijection $\mathfrak{S}_c \to \Sm$,
  denoted $x\mapsto \wt{x}$, is constructed as follows. Starting with $\diag (x)$,
  obtain $\diag (\wt{x})$ by considering only rows indexed by
  $\supp (c)$ and columns  
  indexed by $w(\supp (c))$. Proposition \ref{pr:standard} shows that
  this correspondence is a Bruhat order isomorphism. 

We have $\wt{cw} < \wt{w}$. There is some transposition $\wt{x}\in
\Sm$ such that $\wt{cw} \to \wt{x}\wt{cw} \le \wt{w}$. Observe that
$\wt{x}\wt{cw} = \wt{tcw}$ for some transposition $t\in \Sn$ with
$\supp(t)\subseteq \supp (c)$. Thus, $tuw^{-1}$ has one more cycle than
$uw^{-1}$ does (the cycle $c$ of $uw^{-1}$ is ``split'' upon multiplication by
$t$). It follows that $t$ has the desired properties.
\end{proof}
\end{lemma}

For convenience, let us record as a theorem the various equivalent
conditions that have made appearances in this section.

\begin{theorem}\label{th:collection}
Given a permutation $w \in \Sn$, the following assertions are equivalent:
\begin{itemize}
\item[(i)] $\#\re(w) = \#[e,w]$.
\item[(ii)] $w$ has the right hull property.
\item[(iii)] $w$ avoids the patterns $4231$, $35142$, $42513$ and $351624$.
\item[(iv)] $\al(u,w)=\lp(uw^{-1})$ for all $u\le w$.
\end{itemize}
\begin{proof}
Theorem \ref{th:main} shows (i) $\Leftrightarrow$ (iv), the
  equivalence (ii) $\Leftrightarrow$ (iii) is Sj\"ostrand's Theorem
  \ref{th:sjostrand}, (ii) $\Rightarrow$ (iv) is Lemma 
  \ref{le:RHP} and, finally, (i) $\Rightarrow$ (iii) is the less tricky direction of (B); see \cite[Theorem 4.1]{HLSS}.
\end{proof}
\end{theorem}

\begin{remark} 
{\em A fifth equivalent assertion, which has not been used in this section, was
  given by Gasharov and Reiner in \cite{GR}. They showed that $w\in
  \Sn$ satisfies condition (iii) of Theorem \ref{th:collection} if and
  only if the type $A$ Schubert variety indexed by $w$ is ``defined by
  inclusions'' (see \cite{GR} for the definition). Moreover,
  they discovered that these varieties admit a particularly nice
  cohomology presentation. It would be very interesting to understand
  more explicitly how the other equivalent conditions are connected to this
  picture. Regarding the type independent conditions (i) and (iv),
  this could perhaps lead to interesting cohomological information about
  Schubert varieties of other types.}
\end{remark}

\section{Rational smoothness implies surjectivity}\label{se:ratsmooth}
Suppose $W$ is a Weyl group of a semisimple simply connected complex Lie
group $G$. Then, $W$ is a finite reflection group whose elements index
the Schubert varieties in the (complete) flag variety of $G$. A lot of
work has been devoted to 
understanding how singularities of Schubert varieties are reflected by
combinatorial properties of $W$. A good general reference is \cite{BL}. 

Oh, Postnikov and Yoo established in \cite{OPY} that when $W$ is a
symmetric group, a $q$-analogue of the equality $\#\re(w) = \#[e,w]$
holds whenever the corresponding Schubert variety is rationally
smooth. The same property was conjectured for all Weyl groups
$W$. Recently, Oh and Yoo \cite{OY} presented a proof of this
conjecture. 

In this section, we shall see that the $q=1$ case, i.e.\ the actual
identity $\#\re(w) = \#[e,w]$, of Oh and Yoo's
result is a simple consequence of Theorem \ref{th:main}. In the
process, we formulate a new combinatorial criterion (Theorem
\ref{th:rhombi} below) for detecting rational singularities of
Schubert varieties. 

Let $X(w)$ denote the Schubert variety indexed by
$w\in W$. For the purposes of the present paper, the following classical
criterion could be taken as the definition of $X(w)$ being
rationally smooth. 

\begin{theorem}[Carrell-Peterson \cite{CP}]
The variety $X(w)$ is rationally smooth if and only if the
Bruhat graph $\bg(w)$ is regular, i.e.\ has equally many edges
(disregarding directions) incident with each vertex. 
\end{theorem}
 
For instance $\bg(3412)$, depicted in Figure
\ref{fi:bruhat_graph}, is not regular. Hence, $X(3412)$ is not
rationally smooth. 

If $w\in W$ is understood from the context and $u\le w$, let 
\[
\out(u) = \{t\in T\mid tu \le w\}.
\]
Thus, $\out(u)$ can be thought of as the set of edges incident to $u$ in
$\bg(w)$. Define $\deg(u) = \# \out(u)$. Since $\out(w) = \inv(w)$,
$\deg(w)=\ell(w)$. Hence, $\bg(w)$ is regular if and only if it is
$\ell(w)$-regular.  

\begin{definition} Suppose $x,y,z\le w$. We say that $[e,w]$ contains
  the {\em broken rhombus} $(x,y,z)$ if {\rm (i)} $x \ot y\to z$ and
  {\rm (ii)} $x\to v\ot z$ implies $v\not \le w$.  
\end{definition}

Returning to Figure \ref{fi:bruhat_graph}, several broken
rhombi can be found in $[e,3412]$. One is given by $(2314, 1324,
1342)$, another is $(1432, 1234, 2134)$. 

\begin{theorem}\label{th:rhombi}
The Schubert variety $X(w)$ is rationally smooth if and only if
$[e,w]$ contains no broken rhombi.
\begin{proof}
For a fixed reflection $t\in T$, we partition $T\setminus\{t\}$ in the
following way. For $r\in T\setminus\{t\}$, Let 
\[
C_t(r) = f^{-1}\left(\SP(\{\alpha_r,\alpha_t\})\cap \Phi^+ \right),
\]
where $f:T\to \Phi^+$ is the natural 1--1 correspondence $r\mapsto
\alpha_r$ between reflections and positive roots. In other words,
$C_t(r)$ consists of all reflections that 
correspond to roots in the plane spanned by $\alpha_t$ and
$\alpha_r$, and $\langle C_t(r)\rangle $ is a dihedral reflection subgroup of
$W$. Now, $\{C_t(r)\setminus\{t\}\mid r\in T\setminus\{t\}\}$ is a partition of
$T\setminus\{t\}$.  

Any subgroup of $W$ generated by reflections is a Coxeter group in its
own right with a canonically defined set of Coxeter generators
\cite{dyer}. As was mentioned in the proof of Lemma \ref{le:path}, there
is an isomorphism of directed graphs from the subgraph of 
$\bg(W)$ induced by a coset $\langle C_t(r)\rangle u$ to the Bruhat graph
of the dihedral Coxeter group $D\cong \langle C_t(r)\rangle$. The image of
$[e,w]\cap \langle C_t(r)\rangle u$ is a Bruhat order ideal $I$ in $D$. The
special structure of dihedral Bruhat intervals shows that either the
number of elements of odd respectively of even lengths in $I$ are
equal, or they differ by one. Assuming $I$ contains at least two
elements, in the former case $I$ has a unique
maximum and in the latter it has two maximal elements $m_1\neq m_2$ of
the same Coxeter length. In this case, let $x$ and $z$ be the preimages of $m_1$ and $m_2$, respectively, and choose $y\in \langle C_t(r)\rangle u$ such that $x\ot y \to z$. Then, Dyer's \cite[Lemma 3.1]{dyer} shows that $x \to v \ot z$ implies $v\in \langle C_t(r)\rangle u$. Thus, $(x,y,z)$ forms a broken rhombus in $[e,w]$.   

Observe that in the Bruhat graph of a dihedral group, $u$ and $v$ are
adjacent if and only if $\ell(u)$ and $\ell(v)$ have different parity. 

Suppose $tu\to u \le w$, $t\in T$. If $[e,w]$ contains no broken
rhombi, the above considerations show that $|\out(u)\cap C_t(r)| =
|\out(tu)\cap C_t(r)|$ for all $r\in T\setminus\{t\}$. Thus,
$\deg(tu)=\deg(u)$ so that in fact all vertices in $[e,w]$ have degree
$\deg(w)$, and $X(w)$ is rationally smooth by the
Carrell-Peterson criterion.

For the converse statement, assume $(x,y,z)$ is a broken rhombus in
$[e,w]$ with $\ell(y)$ maximal. Let $t=xy^{-1}$ and $r=zy^{-1}$. Then, $y$ has one more neighbour in $\langle C_t(r)\rangle y$ than $x$ does. That is,
$|\out(y)\cap C_t(r)| = |\out(x)\cap C_t(r)| + 1$. Moreover, by
maximality of $y$, there is no $r^\prime\in T$ with
$|\out(y)\cap C_t(r^\prime)| = |\out(x)\cap C_t(r^\prime)| -
1$. Therefore $\deg(y)>\deg(x)$, implying that $X(w)$ is rationally
singular.
\end{proof}
\end{theorem}


With this criterion and Theorem \ref{th:main} at our disposal, the $q=1$ case of Oh and Yoo's result is little more than an observation:

\begin{corollary}
The map $\phi$ is surjective, hence bijective, if $X(w)$ is rationally smooth.
\begin{proof}
Suppose $\phi$ is not surjective. Assume $z\le w$ is such that
$\al(z,w)>\lp(zw^{-1})$ and $\lp(zw^{-1})$ is minimal among all
$z$ with this property. By Lemma \ref{le:path}, there exist
$x,y\le w$ such that $x\ot y\to z$ and $\lp(xw^{-1}) = \lp(yw^{-1})-1
= \lp(zw^{-1})-2$. Now, $x\to v\ot z$ implies $v\not \le w$; otherwise
a directed path of length $\al(v,w)+1 = \lp(vw^{-1})+1 \le \lp(xw^{-1})+2$ would exist
from $z$ to $w$, contradicting the assumptions. Hence, $(x,y,z)$ is a
broken rhombus. Theorem \ref{th:rhombi} concludes the proof.  
\end{proof}
\end{corollary}

\end{document}